\documentclass[12pt]{article}
\usepackage{fullpage}
\begin{document}
\newcommand{\qed}{\hbox to 0pt{}\hfill$\rlap{$\sqcap$}\sqcup$ \vspace{3mm}}
\def\NN{\hbox{I\kern-.2em\hbox{N}}}
\def\RR{\hbox{I\kern-.2em\hbox{R}}}
\def\CC{\hbox{I\kern-.5em\hbox{C}}}

\title{New Stability Conditions for Linear   
Difference Equations using Bohl-Perron 
Type Theorems}
\author{Leonid Berezansky $^{1}$   \\
Department of Mathematics,
Ben-Gurion University of the Negev, \\
Beer-Sheva 84105, Israel
\\
and  Elena Braverman $^{2}$ \\
Department of Mathematics and Statistics, University of Calgary, \\
2500 University Drive N.W., Calgary, AB T2N 1N4, Canada}
\date{}

\maketitle

\footnotetext[1]{Partially supported by Israeli Ministry of Absorption}
\footnotetext[2]{Partially supported by the NSERC Research Grant}


\medskip



\begin{abstract}
The Bohl-Perron result on exponential dichotomy 
for a linear difference equation 
$$
x(n+1)-x(n) + \sum_{l=1}^m a_l(n)x(h_l(n))=0, ~~ h_l(n)\leq n,
$$
states (under some natural conditions) that if
all solutions of the non-homogeneous equation with a bounded
right hand side are bounded, then the relevant homogeneous equation
is exponentially stable. According to its corollary, if a given
equation is {\em close} to an exponentially stable comparison equation
(the norm of some operator is less than one), then the considered equation
is exponentially stable.  

For a difference equation with several variable delays and coefficients 
we obtain new exponential stability tests using the above results, 
representation of solutions and comparison equations with a positive 
fundamental function. 
\end{abstract}

\newtheorem{uess}{Lemma}
\newtheorem{guess}{Theorem}
\newtheorem{corol}{Corollary}
\newtheorem{example}{Example}

\section{Introduction}

In this paper we study stability of a scalar
linear difference equation with several delays
\begin{equation}
\label{1}
x(n+1)-x(n)=-\sum_{l=1}^m a_l(n)x(h_l(n)), ~~ h_l(n)\leq n,
\end{equation}
where $h_l(n)$ is an integer for any $l=1, \cdots, m$ and $n=0,1,2, \cdots 
~$
under the following two restrictions on the parameters of 
(\ref{1}) which mean that coefficients and delays are bounded:

{\bf (a1)}  there exists $K>0$ such that ~~~ ${\displaystyle 
 |a_l(n)| \leq K }$ for $l=1, \cdots, m$, $n=0,1,2, \cdots$ ;

{\bf(a2)} ~there exists $T>0$ such that $n-T \leq h_l(n) \leq n$
for $l=1, \cdots, m$, $n=0,1,2, \cdots$ ~~.
\vspace{1mm}

Stability of equation (\ref{1}) has been
an intensively developed area during the last two decades, see
\cite{JMAA2005,BBL}, \cite{n20}-\cite{GyoriHartung}, \cite{GLV,n21},
\cite{KP}-\cite{ZhTianChuan} 
and references therein.

In the present paper we study a connection between stability and 
existence of a positive solution for a general linear difference equation 
with nonnegative variable coefficients and several delays.
This idea was developed in \cite{GLV} for equations with constant 
coefficients, see also \cite{Pituk2003} and references therein for
some further results on equations with variable coefficients and
a nonlinear part.

The method of the present paper is based on the Bohl-Perron type theorems
which connect the boundedness of all solutions for all bounded right hand 
sides with the exponential stability of the relevant homogeneous equation.
Here we apply the development of this method \cite{FDE2004,JMAA2005} where
stability properties of the original equation are established based on
the known asymptotics of an auxiliary (comparison, model) equation
(see Lemma \ref{basiclemma} below). This idea 
can be compared to \cite{Pituk2003} where a nonlinear perturbation of a linear
equation is considered. The difference between the present paper and
\cite{JMAA2005} is that the model equation is, generally, a high order 
equation, which allows to consider equations with large delays. 

The paper is organized as follows. In Section 2 we deduce some general 
exponential stability results for high order difference equations with 
variable coefficients. Section 3 presents explicit stability tests in
terms of delays and coefficients. Finally, Section 4 involves  
discussion and examples, which compare our results with known 
stability tests, and outlines open problems.

\section{Exponential Stability: General Results}

We assume that (a1)-(a2) hold for equation (\ref{1}) 
and all other equations in the paper.
In particular, the system has a finite prehistory: 
$h_l(n) \geq n_0 - T$, $n \geq n_0$, for any $n_0 \geq 0$.
We note that, generally, under this assumption (\ref{1}) can be
written as a higher order equation ${\displaystyle x(n+1)-x(n)=
-\sum_{l=0}^T b_l(n)x(n-l) }$, where $b_l(n)$ are defined by coefficients 
$a_k(n)$ of equation (\ref{1}). 
However, there are several reasons 
why we prefer form (\ref{1}).
\begin{enumerate}
\item
The number of nonzero terms in (\ref{1}) can be much smaller than
$T$; this form allows us to specify stability conditions for equations
with one, two or three variable delays, equations with positive and 
negative coefficients (where positive and negative terms are written 
separately) and so on. In particular, we refer to the long tradition
when, for example, the equation with one variable delay 
${\displaystyle x(n+1)-x(n)= a(n)x(h(n))}$ was studied in this form.
\item
The main idea of the method applied in the present paper 
is based on the comparison to a model stable equation.
Form (\ref{1}) gives some flexibility: in many cases, the
model equation may have different delays compared to the original
equation.  
\end{enumerate}

In future, we will consider the scalar linear difference 
equation
\begin{equation}
\label{1gen}
x(n+1)-x(n)=-\sum_{l=1}^m  a_l(n)x(h_l(n)) +f(n),~~ n \geq n_0,
\end{equation}
\begin{equation}
\label{2}
x(n)=\varphi(n),~~ n\leq n_0,
\end{equation}
and the corresponding homogeneous equation (\ref{1}).

Further we will extensively apply
the solution representation formula and 
properties of the fundamental function.
We start with the definition of this function.
\vspace{2mm}

\noindent
{\bf Definition}.
The solution $X(n,k)$ of the problem
$$
x(n+1)-x(n)=-\sum_{l=1}^m a_l(n)x(h_l(n)),~ n>k, 
~~~x(n)=0,~n<k,~~x(k)=1
$$
is called {\em the fundamental function} of equation (\ref{1gen})
(or of ~(\ref{1}) ).

In the following we will use a modification of the solution 
representation formula \cite{E}
\begin{equation}
\label{3add}
x(n)=X(n,n_0)x(n_0)+\sum_{k=n_0}^{n-1} X(n, k+1)f(k)-\sum_{k=n_0}^{n-1}
X(n,k+1)\sum_{l=1}^m
a_l(k)\varphi(h_l(k)),
\end{equation}
where $\varphi(h_l(k))=0,h_l(k)\geq n_0$.



Denote by ${\bf l}^{\infty}$  the space of bounded
sequences $v=\{ x(n) \}$  with the norm  \\ ${\displaystyle  
\| v \|_{{\bf l}^{\infty}}=\sup_{n \geq 0} |x(n)|<\infty }$.

\begin{uess} $\cite{JMAA2005}$
\label{basiclemma}
Suppose for the equation ${\cal L}(\{ x(n) \})=0$, where
$$
{\cal L}\left( \{ x(n) \}_{n=0}^{\infty} \right) = \left\{ 
x(n+1)-x(n) + \sum_{l=1}^m a_l(n)x(h_l(n)) \right\}_{n=0}^{\infty},
$$
conditions (a1)-(a2) hold, 
the equation ${\cal L}_1(\{ x(n) \})=0$ is exponentially stable,  where
\begin{equation}
\label{comparison}
{\cal L}_1\left(\{ x(n) \}_{n=0}^{\infty} \right) = \left\{ x(n+1)-x(n) + 
\sum_{l=1}^r b_l(n)x(g_l(n)) \right\}_{n=0}^{\infty}.
\end{equation} 
Let $Y(n,k)$ be the fundamental function of 
(\ref{comparison}) and ${\cal C}_1$ be {\em the Cauchy operator}
of equation (\ref{comparison}) which is
\begin{equation}
\label{Cauchy}
{\cal C}_1 \left\{ f(n) \right\}_{n=0}^{\infty} = \left\{ y(n)= 
\sum_{k=0}^{n-1} Y(n,k+1)f(k) \right\}_{n=0}^{\infty},
\end{equation}
with $y(0)=0$. If the operator ${\cal T} = {\cal C}_1 ({\cal L}_1-{\cal 
L}) = I - {\cal C}_1 {\cal L}$ satisfies $\| {\cal T} \|_{{\bf l}^{\infty}
\to {\bf l}^{\infty}} < 1$, then the equation
${\cal L}(\{ x(n) \})=0$ is exponentially stable.
\end{uess}

\noindent
{\bf Definition.} Equation (\ref{1}) is {\em exponentially 
stable} if there exist constants $M>0$,
$\lambda\in (0,1)$ such that for every solution $\{x(n)\}$ of
(\ref{1}),(\ref{2}) the inequality
\begin{equation}
\label{17add}
|x(n)|\leq M \lambda^{n-n_0} \left(\max_{ -T+n_0 \leq k \leq n_0}
\left\{|\varphi(k)|\right\}\right)
\end{equation}
holds for all $n\geq n_0$, 
where $M, \lambda$ do not depend on $n_0$.
\vspace{2mm}



The following result states that the exponential estimate of the 
fundamental function implies exponential stability (see 
\cite{JMAA2005}, Theorem 4).

\begin{uess} $\cite{JMAA2005}$
\label{lemma1b}
Suppose there exist 
$L, \mu$, $L>0$, $0<\mu < 1$, such that
\begin{equation}
\label{star}
|X(n,k)| \leq L \mu^{n-k}, ~~ n \geq k \geq n_0.
\end{equation}
Then there exist $M, \lambda$, $M>0$, 
$0<\lambda <1$
such that inequality 
(\ref{17add}) holds for a solution of (\ref{1}), 
(\ref{2}).
Conversely, if (\ref{17add}) holds, then (\ref{star}) is also valid, 
with $L=M$, $\mu=\lambda$.
\end{uess}

We proceed to estimation of the fundamental function.
Let us note that throughout the paper we assume that the sum equals zero
and the product equals one if the lower bound exceeds the upper bound.

\begin{uess}
\label{est}
The fundamental function of equation (\ref{1}) has the following estimate
$$
|X(n,k)| \leq \prod_{j=k}^{n-1} \left( 1+ \sum_{l=1}^m \left| 
a_l(j)\right|
\right).
$$
\end{uess}
{\bf Proof.}
Let us fix $k$ and denote $x(n)=X(n,k)$. By definition
we have
$$
x(n)=x(n-1)- \sum_{l=1}^m a_l(n-1) x(h_l(n-1)),
$$ $$
x(n)=0, ~n<k,~~ x(k)=1.
$$
Then
$$
|x(n)| \leq |x(n-1)|+ \sum_{l=1}^m \left| a_l(n-1)\right|
\left| x(h_l(n-1)) \right|.
$$
Consider the linear equation
$$
y(n)=y(n-1)+ \sum_{l=1}^m \left| a_l(n-1) \right| y(h_l(n-1))
$$
with the initial condition
$$  
y(n)=0, ~n<k,~~ y(k)=1.
$$
By induction it is easy to see that $\{y(n)\}$ is a positive nondecreasing
sequence and
$|x(n)|\leq y(n)$, so $y(h_l(n-1)) \leq y(n-1)$.
We have
$$
y(n)\leq y(n-1)\left( 1+ \sum_{l=1}^m \left| a_l(n-1) \right| \right)
\leq
\prod_{j=k}^{n-1} \left( 1+ \sum_{l=1}^m  \left| a_l(j) \right| 
\right)y(k)
=\prod_{j=k}^{n-1} \left( 1+ \sum_{l=1}^m  \left| a_l(j) \right| \right). 
$$
The remark that
$|X(n,k)| = |x(n)|\leq y(n)$ completes the proof.
\qed

The following Lemma claims that under (a1)-(a2) exponential stability is 
really an asymptotic property in the sense that if $a_l(k)$ or $X(n,k)$ 
is changed on any finite segment $n_0\leq k \leq n_1$, this does not 
influence stability properties of (\ref{1}).

\begin{uess}
\label{lemma1a}
If for  some $n_1 > n_0$ there exist $L_1$, $\mu$, $L_1>0$,
$0<\mu <1$, such that
\begin{equation}
\label{2star}   
|X(n,k)| \leq L_1 \mu^{n-k}, ~~ n \geq k \geq n_1,
\end{equation}
then for some $L$, $L>0$, inequality (\ref{star}) is satisfied.
\end{uess}
{\bf Proof.} 
Denote
\begin{equation}
\label{denotA}  
A=\mu^{n_0-n_1} \prod_{j=n_0}^{n_1} \left( 1+
\sum_{l=1}^m |a_l(j)| \right).
\end{equation} 
For $n_0 \leq k \leq n \leq n_1$ we have $n-k \leq n_1-n_0$ and
$0<\mu<1$, so Lemma \ref{est} implies
$$
|X(n,k)| \leq \prod_{j=k}^{n-1} \left( 1+ \sum_{l=1}^m |a_l(j)| \right)
\leq A \mu^{n_1-n_0}\leq A \mu^{n-k},
~~n_0\leq k \leq n \leq n_1.$$

If $n_0\leq k\leq n_1$, $n>n_1$, then by (\ref{3add}), (\ref{2star}), (a1) 
and (a2)
we have $h_l(j)\geq n_1$, $j>n_1+T$, thus (the second sum below involves 
only such $X(h_l(j),k)$ that $h_l(j)<n_1$)
$$|X(n,k)| \leq |X(n,n_1)| |X(n_1,k)|+ \sum_{j=n_1}^{n_1+T} |X(n,j+1)|
\sum_{l=1}^m |a_l(j)| |X(h_l(j),k)| $$
$$ \leq L_1 \mu^{n-n_1}A\mu^{n_1-k} + \sum_{j=n_1}^{n_1+T} L_1 
\mu^{n-(j+1)}
\left(\sum_{l=1}^m |a_l(j)| A\right)
\leq AL_1 \mu^{n-k}+ K m A L_1 \sum_{j=n_1}^{n_1+T}
\mu^{n-(j+1)} $$
$$ =AL_1\mu^{n-k}+K m A L_1 \mu^{n-n_1}\sum_{j=0}^{T}\mu^{-(j+1)}
\leq B \mu^{n-k}, $$
where  ${\displaystyle B= AL_1 + K m A L_1 \mu^{n_0-n_1}\sum_{j=0}^{T}
\mu^{-(j+1)} }$ does not depend on $n$, $A$ was
defined in (\ref{denotA}). Choosing $L = \max \{ A,B \}$,
we obtain estimate (\ref{star}) for $n \geq k \geq n_0$.
\qed

Consider the  equation
\begin{equation}
\label{dstar}
x(n+1)-x(n) = - \sum_{l=1}^m b_l(n) x(h_l(n)), ~~n \geq n_0.
\end{equation}

\begin{corol}
\label{corol1}
Suppose $a_l(n)= b_l(n),~~  n\geq n_1> n_0$.
Equation (\ref{dstar}) is exponentially stable
if and only if (\ref{1}) is exponentially stable.
\end{corol}
{\bf Proof.}
Let $X(n,k)$ and $Y(n,k)$ be  fundamental functions of (\ref{1}), 
and (\ref{dstar}), respectively. Suppose (\ref{1}) is exponentially 
stable.
Then for $X(t,s)$ inequality (\ref{star}) holds with $L>0, 0<\mu<1$.
Since $X(n,k)=Y(n,k)$ for $n\geq k\geq n_1$, then $|Y(n,k)| \leq L
\mu^{n-k}$, $n \geq k \geq n_1.$ 

By Lemma \ref{lemma1a} we have $|Y(n,k)| \leq L_1
\mu^{n-k}, ~~ n \geq k
\geq n_0$, for some $L_1 >0$.
Finally, Lemma \ref{lemma1b} implies that equation(\ref{dstar}) is 
exponentially stable.

The same argument proves the converse statement.
\qed

Now we proceed to a result on exponential stability of
equation (\ref{1}) with a positive fundamental
function and nonnegative coefficients.

\begin{guess}
\label{theorem3}
Suppose the fundamental function of
(\ref{1}) is eventually positive, i.e., for some $n_0 \geq 0$ we have
$X(n,k)>0,~ n \geq k \geq n_0$, and, in addition, either
\begin{equation}
\label{16add}
a=\liminf_{n \to \infty} \sum_{l=1}^m a_l(n) >0
\end{equation}
or a more general condition holds: \\ there exists a positive integer $p$, 
such that
\begin{equation}
\label{16a}
b= \limsup_{n \to \infty} \prod_{j=n}^{n+p-1}
\left( 1 - \sum_{l=1}^m a_l(j) \right) <1.
\end{equation}

Then each of the following statements is valid.
\begin{enumerate}
\item
The fundamental function of (\ref{1}) is eventually nonincreasing.
\item
If (\ref{16add}) holds, then  (\ref{star}) is satisfied for any
$\mu$, $1-a< \mu <1$;
if (\ref{16a}) holds, then exponential estimate (\ref{star}) is valid for
any $\mu, ~b^{1/p}<\mu<1$.
\item
Equation (\ref{1}) is exponentially stable.
\end{enumerate} 
\end{guess}
{\bf Proof.}
By the assumptions of the theorem for $n \geq n_0+T$, $k \geq n_0$
we have
$$ X(n+1,k)=X(n,k) - \sum_{l=1}^m a_l(n) X(h_l(n), k) \leq X(n,k),$$
since $a_l(n)$, $X(h_l(n), k)$ are nonnegative, so $X(n,k)$ is
nonincreasing for $n \geq k \geq n_0+T$, which completes the proof of the
first part of the theorem.

Further, let us assume that (\ref{16add}) holds, then there exist
$\varepsilon$, $0< \varepsilon<a$, and $n_1 \geq n_0+T$ such that
$$
\sum_{l=1}^m a_l(n) \geq \varepsilon, ~~n \geq n_1.
$$

Thus $X(n,k)$ is nonincreasing in $n$ for any
fixed $k$, $n \geq k \geq n_1$. Consequently, for $n\geq k+T$ we have
$X(h_l(n),k) \geq X(n,k)$ and
$$X(n+1,k) \leq X(n,k)- \sum_{l=1}^m a_l(n) X(n, k) =
\left( 1 - \sum_{l=1}^m a_l(n) \right) X(n,k).$$
The latter inequality yields that ${\displaystyle 0\leq \sum_{l=1}^m 
a_l(n)<1}$, $n \geq n_2=n_1+T$, which implies $a \leq 1$.
Thus $X(k+1,k)  <  (1-\varepsilon) X(k,k) = 1-\varepsilon$. Repeating
this procedure, we obtain 
$$X(n,k) \leq \mu^{n-k},~~n \geq k \geq n_1,
$$
where $\mu=1-\varepsilon$. By Lemma \ref{lemma1a} estimate (\ref{star}) is
valid for $n\geq k \geq n_0$, with the same $\mu$ and some $L>0$.

Next, let (\ref{16a}) hold and the expression under $\limsup$ be less than
$q<1$ for $n \geq n_1$, where $n_1 \geq n_0+T$.
Then the corresponding estimate for any positive integer is
$X(k+rp,k) \leq q^r,$
so $$X(n,k) \leq L \mu^{n-k}, ~~n \geq k \geq n_1,~  
\mbox{ where ~} \mu=q^{1/p}, ~~L=\mu^{-p}.$$
Lemma \ref{lemma1a} implies the
exponential estimate for the fundamental function, $n \geq k \geq 
n_0$, thus the proof of the second part is complete.

Finally, by Lemma \ref{lemma1b} equation (\ref{1}) is exponentially
stable.
\qed

\noindent
{\bf Remark 1.}
Let us remark that $a\leq 1$ in (\ref{16add}) since all sums are 
less than one,
beginning with some $n_2=n_1+T$. If $a=1$, then an exponential estimate 
with any positive $\mu$ will work, which is illustrated by the following
example.

\begin{example} 
For the equation
$$x(n+1)-x(n)=- \left( 1-\frac{1}{n+1} \right) x(n), ~~n \geq 0,$$  
we have $a=1$, where $a$ is defined in (\ref{16add}), the fundamental
function  is $X(n,k)=k!/n!$\,. For this equation
exponential estimate (\ref{star}) is valid for any $\mu$,
$0<\mu<1$, since
$$\frac{k!}{n!} \leq \mu^{n-k}, ~~ n \geq k \geq \frac{1}{\mu}.$$
\end{example}

\noindent
{\bf Remark 2.} Let us notice that under the assumptions of Theorem 
\ref{theorem3} inequality (\ref{16add}) implies (\ref{16a}) 
for any positive integer $p$ and for $b$ in (\ref{16a})
which does not exceed $(1-a)^p<1$
(here we assume that all other conditions of the theorem are satisfied).
Another inequality which implies
(\ref{16a}) is
the existence of a positive integer $p$ such that
\begin{equation}
\label{16b}
c= \liminf_{n \to \infty} \max_{n \leq k \leq n+p-1} \sum_{l=1}^m a_l(k)
>0.
\end{equation}
If (\ref{16b}) holds then for any $\varepsilon>0$ and $n$ large enough
among $p$ successive sets of coefficients at least one satisfies 
${\displaystyle \sum_{l=1}^m a_l(n)> \frac{c}{p}-\varepsilon }$, so 
(\ref{16a})
is valid with $b \leq  1-c/p<1$.

Thus, later we will refer to condition (\ref{16a}) only.

Further, we will apply nonoscillation tests, the following 
result is Theorems 4.1 in \cite{ADSA}.

\begin{uess} $\cite{ADSA}$
\label{ADSA51}
Suppose $a_l(n)\geq 0$, $l=1,2, \cdots , m $, and
for some $n_0\geq 0$
\begin{equation}
\label{10}
\sup_{n\geq n_0} \sum_{l=1}^m
a_l (n)< \frac{1}{2}, ~~~
\sup_{n\geq n_0}\sum_{l=1}^m ~~ \sum_{k=\max\{n_0,\min_l h_l(n)\}}^{n-1}
a_l(k)\leq \frac{1}{4}.
\end{equation}
Then the fundamental function of (\ref{1}) is eventually positive: 
$X(n,k)>0, ~n\geq n_0$.
\end{uess}

We remark (\cite{GL}, Theorem 7.2.1) that the condition
\begin{equation}
\label{autonomnonoscil}
a \leq \frac{k^k}{(k+1)^{(k+1)}}
\end{equation}
is necessary and sufficient for nonoscillation of the autonomous equation
\begin{equation}
\label{autonomous}
x(n+1)-x(n)=-ax(n-k), ~k\geq 1.
\end{equation}

\begin{corol}
\label{corol2}
Suppose $a_l(n) \geq 0$, (\ref{16a}) holds and (\ref{10}) is satisfied for
some $n_0\geq 0$.
Then (\ref{1}) is exponentially stable.
\end{corol}
{\bf Proof.}
By Lemma \ref{ADSA51} inequalities (\ref{10})
imply that the fundamental function
of (\ref{1}) is eventually positive.
Application of Theorem \ref{theorem3} completes the proof.
\qed

Consider together with (\ref{1}) the following comparison 
equation
\begin{equation}
\label{7add}  
x(n+1)-x(n)=-\sum_{l=1}^m b_l(n)x(g_l(n)), ~n \geq  n_0,
\end{equation}
where $g_l(n)\leq n$. Denote by $Y(n,k)$ the fundamental function of
equation (\ref{7add}).

\begin{uess} $\cite{ADSA}$
\label{ADSA41}
Suppose $a_l(n)\geq b_l(n)\geq 0,
~g_l(n)\geq h_l(n)$, $l=1,2, \cdots , m$,
for sufficiently large $n$.
If equation (\ref{1}) has an eventually positive solution, then 
(\ref{7add})
has an eventually
positive solution and its fundamental function $Y(n,k)$ is eventually 
positive.
\end{uess}

\begin{corol}
\label{corol3}
Suppose (\ref{16a})   and at least one of the following conditions hold:

1) $ 0 \leq a_l(n) \leq \alpha_l$, $h_l(n)\geq n-\tau_l$
and there exists $\lambda>0$ such that
\begin{equation}
\label{20}   
\lambda-1\leq -\sum_{l=1}^m \alpha_l \lambda^{-\tau_l}.
\end{equation}   

2) $m=1, n-h(n)\leq k, a(n)\leq \frac{k^k}{(k+1)^{(k+1)}}$.

Then (\ref{1}) is exponentially stable.
\end{corol}
{\bf Proof.}
Suppose $\lambda_0>0$ is a solution of inequality (\ref{20}) and denote
${\displaystyle f(\lambda)=\lambda-1+\sum_{l=1}^m \alpha_l
\lambda^{-\tau_l} }$.
We have $f(\lambda_0)\leq 0,~ f(1)\geq 0$, hence there exists
a positive solution
$\lambda_1>0$ of the characteristic equation for
the autonomous equation
$$
x(n+1)-x(n)=-\sum_{l=1}^m \alpha_l x(n-\tau_l),
$$
which consequently has  a positive solution $x(n)=\lambda_1^n$.
Lemma \ref{ADSA41} implies that (\ref{1}) also has a
positive solution. By Theorem \ref{theorem3} equation (\ref{1}) is
exponentially stable.

Proof of the second part is similar.
\qed

\noindent   
{\bf Remark 3.} By Theorem 3.1 in \cite{ADSA}
it is enough to assume the existence of an eventually positive solution
in the conditions of Theorem \ref{theorem3} 
rather than to require that the fundamental function is positive.

\section{Explicit Stability Tests} 

%

Lemma \ref{basiclemma} claims that if 
(\ref{1}) is in some sense close to an exponentially stable equation, 
then it is also exponentially stable. Lemma \ref{fundamental_1} provides 
some estimates which are useful to establish this closeness. 

Further, we deduce explicit exponential 
stability conditions based on Lemmas \ref{basiclemma} and 
\ref{fundamental_1}.  As above, we assume that (a1)-(a2) hold for 
(\ref{1}).

\begin{uess}
\label{fundamental_1}
Suppose the fundamental function of the equation 
(\ref{1}) is positive:
$X(n,k)>0$, $n \geq k \geq n_0$, and $a_l(n) \geq 0$, $l=1, \cdots , m$,
$n \geq n_0$.
Then there exists $n_1 \geq n_0$ such that
\begin{equation}
\label{51}
0 \leq \sum_{k=n_0}^{n-1} X(n,k+1) \sum_{l=1}^m a_l(k) \leq 1, ~~~
n \geq n_1.
\end{equation}
\end{uess}
{\bf Proof.}
Since $a_l(n) \geq 0$ and $X(n,k)>0$ for $n \geq k 
\geq n_0$, then $X(n,k)$ is nonincreasing in $n$ for any $k$,
$n \geq k \geq n_0$. Consider
\begin{equation}
\label{52}
z(n)= \left\{ \begin{array}{ll} 1, & n \geq n_0, \\ 0, & n<n_0.
\end{array}  \right.
\end{equation}
Then 
$$z(n+1)-z(n)+ \sum_{l=1}^m a_l(n) z(h_l(n))=f(n), ~~ n \geq n_0,$$ with
${\displaystyle
f(n)= \sum_{l=1}^m a_l(n) \chi_{n_0} (h_k(n))}$, where
${\displaystyle \chi_{n} (j)= \left\{ \begin{array}{ll} 1, & j \geq n, \\ 
0, & j<n.  
\end{array} \right. }$. 

Thus by the solution representation (\ref{3add}) for $n>n_0$ we have
$$1=z(n)= X(n,n_0)+\sum_{k=n_0}^{n-1} X(n,k+1)f(k) =
X(n,n_0)+\sum_{k=n_0}^{n-1} X(n,k+1) \sum_{l=1}^m a_l(k) \chi_{n_0} 
(h_l(k)).
$$
Since by (a2) all delays are bounded, then there exists a maximal 
delay 
$T$ such that $\chi_{n_0} (h_l(k))=1$ for any $l=1, \cdots, m$ and $k\geq 
n_0+T$. Thus for $k\geq n_0+T$ we have
$$0 \leq \sum_{k=n_0}^{n-1} X(n,k+1)  \sum_{l=1}^n a_l(k) = 1- 
X(n,n_0)<1.$$
Since $X(n,k)$ is nonincreasing and positive, then
$0<X(n,n_0) \leq 1$, $n >n_0$. Thus ${\displaystyle 
0 \leq \sum_{k=n_0}^{n-1} X(n,k+1) \sum_{l=1}^m a_l(k) \leq 1, ~n \geq 
n_1=n_0+T}$, which completes the proof.
\qed

Now let us proceed to explicit stability conditions.

\begin{guess}
\label{dominate}
Suppose there exists a subset of indices
$I \subset \{ 1,2, \cdots , m\}$ such that $a_k \geq 0$, $k \in I$,
for the sum 
${\displaystyle \sum_{l \in I} a_l(n) }$ inequality (\ref{16a}) holds,
the fundamental function $X_1(n,k)$ of the equation
\begin{equation}
\label{54}
x(n+1)-x(n)+ \sum_{l \in I} a_l(n) x(h_l(n))=0
\end{equation}
is eventually positive and 
\begin{equation}
\label{55}
\limsup_{n\to \infty} \frac{\sum_{l \not\in I} |a_l(n)| }{
\sum_{l \in I} a_l(n)} < 1.
\end{equation}
Then equation (\ref{1}) is exponentially stable.
\end{guess}
{\bf Proof.} We apply the same method as in \cite{JMAA2005}, where the 
comparison (model) equation is (\ref{54}).
\qed

\noindent
{\bf Remark 4.} 
Most known explicit stability results include estimates where
coefficients $a_l(k)$  are summed up in $k$ from $h_l(n)$ to $n$.
We note that if the comparison equation is
 \begin{equation}
\label{18}
x(n+1)-x(n)=-b(n)x(n-k)
\end{equation}
then the upper index is $n-k-1$.

Now we will take general exponentially stable difference
equations with a positive fundamental function as a class of comparison
equations.
\begin{corol}
\label{theorem5} 
Suppose there exist a set of indices $I \subset \{ 1,2,
\cdots, m \}$, functions $g_l(n)\leq n, l\in I$, and
positive numbers $\alpha_0$, $\alpha_1$ and
$\gamma<1$, such that for $n$ 
sufficiently large the inequalities
$$
0< \alpha_0 \leq  \sum_{l \in I} a_l(n)\leq \alpha_1 < 1,~ l\in I,
$$
hold and the difference equation
\begin{equation}
\label{21}
x(n+1)-x(n)=-\sum_{l \in I} a_l(n)x(g_l(n))
\end{equation}
has a positive fundamental function. If
\begin{equation}
\label{22}
\sum_{k \in I} | a_k(n) |
\sum_{j=\min\{h_k(n),g_k(n)\}}^{\max\{h_k(n),g_k(n)\}-1} \sum_{l=1}^m
|a_l(j)| + \sum_{k \not\in  I} |a_k(n)| \leq \gamma \sum_{k \in I}  
a_k(n)
\end{equation}
then (\ref{1}) is exponentially stable.
\end{corol}
{\bf Proof.}
Let us rewrite (\ref{1}) in the form
$$x(n+1)-x(n)=-\sum_{k \in I} a_k(n) x(g_k(n))
+ \sum_{k \in I} a_k(n) [x(g_k(n))-x(h_k(n))]  - \sum_{k \not\in  I}  a_k(n) 
x(h_k(n))
$$ $$=
-\sum_{k \in I} a_k(n) x(g_k(n))
+ \sum_{k \in I} a_k(n) \sigma_k 
\sum_{j=\min\{h_k(n),g_k(n)\}}^{\max\{h_k(n),g_k(n)\}-1} 
[x(j+1)-x(j)] - \sum_{k \not\in  I}  a_k(n)
x(h_k(n))
$$
$$ =
-\sum_{k \in I} a_k(n) x(g_k(n)) - \sum_{k \in I} a_k(n) \sigma_k 
\sum_{j=\min\{h_k(n),g_k(n)\}}^{\max\{h_k(n),g_k(n)\}-1}
\sum_{l=1}^m a_l(j)x(h_l(j))-\sum_{k \not\in  I}  a_k(n) x(h_k(n)),
$$
where ${\displaystyle \sigma_k=\left\{ \begin{array}{ll} 1, 
& g_k(n)>h_k(n), \\ -1 & g_k(n)<h_k(n) \end{array}\right. }$.
Since (\ref{21}) has a positive fundamental function, then the reference
to Theorem \ref{dominate} completes the proof.
\qed

\noindent
{\bf Remark 5.} Based on the choice of subset $I$,
the theorem contains $2^m-1$ different stability 
conditions.
\vspace{2mm}

Assuming 
$I=\{ 1,2, \cdots , m\}$
in Corollary \ref{theorem5}, we obtain the following result.

\begin{corol}
\label{corol4}
Suppose there exists $g(n)\leq n$ and
positive numbers $\alpha_0$, $\alpha_1$,  $\gamma<1$ such that for $n$ 
sufficiently large
\begin{equation}
\label{25}
0< \alpha_0 \leq b(n):=\sum_{l =1}^m a_l(n)\leq \alpha_1<1
\end{equation}
and the difference equation
\begin{equation}
\label{26}
x(n+1)-x(n)=-\sum_{l=1}^m a_l(n)x(g(n))
\end{equation}   
has a positive fundamental function. 
If  for $n$ large enough
$$
\sum_{l=1}^m |a_l(n) | 
\sum_{k=\min\{h_l(n),g(n)\}}^{\max\{h_l(n),g(n)\}-1} \sum_{l=1}^m |a_l(k)|
\leq \gamma \sum_{l=1}^m a_l(n), $$
then equation (\ref{1}) is exponentially stable.
\end{corol}

The following result is an immediate corollary 
of Lemma \ref{ADSA51} and Theorem \ref{dominate}.

\begin{corol}  
\label{theorem2}
Suppose ${\displaystyle 0< a_0 \leq a_0(n) \leq 
b_0<\frac{1}{4} }$ and there exists $\gamma$ such that $0<\gamma<1$ and 
${\displaystyle \sum_{l=1}^m |a_l(n) | \leq \gamma 
a_0(n)}$
for $n$ large enough. Then the equation
\begin{equation}
\label{19}
x(n+1) - x(n) = -a_0(n)x(n-1) - \sum_{l=1}^m a_l(n)x(h_l(n))
\end{equation}
is exponentially stable.
\end{corol}


\begin{corol}
\label{theorem3add}
Suppose for some positive $a_0, b_0, \gamma$,
where $b_0<1, \gamma<1$, the following inequalities are satisfied for $n$
large enough
\begin{equation}
\label{add1}
0< a_0 \leq  \sum_{l=1}^m a_l(n)  \leq b_0<1/4,
\end{equation}
$$
\sum_{k=1}^m |a_k(n) | \sum_{j=h_{k}(n)}^{n-2} 
\sum_{l=1}^m \left| a_k(j) \right|
\leq \gamma \sum_{l=1}^m a_l(n). $$
Then equation (\ref{1}) is exponentially stable.
\end{corol} 


Now let us consider the case $m=2$
\begin{equation} 
\label{27}
x(n+1)-x(n) = -a(n)x(g(n)) - b(n)x(h(n)).
\end{equation}

\begin{corol} 
\label{corol5}
Suppose 
there exist $a_0>0$ and $\gamma$,
$0<\gamma<1$ such that at least one of the following conditions holds 
for $n$ sufficiently large:\\
1) ${\displaystyle
 0<a_0 < a(n)\leq a_1<\frac{1}{2},~  \sum_{k= g(n)}^{n-1}
a(k)\leq \frac{1}{4}, ~ |b(n)| \leq \gamma a(n);
}$ \\ \\
2) ${\displaystyle
0<a_0 \leq a(n)+b(n) \leq a_1< \frac{1}{2}, ~\sum_{k=g(n)}^{n-1}
(a(k)+b(k))\leq \frac{1}{4},}$ 
$$ 
|a(n)| \sum_{k=\min\{h(n),g(n)\}}^{\max\{h(n),g(n)\}-1}
[|a(k)|+|b(k)|] < \gamma [a(n)+b(n)].
$$
Then equation (\ref{27}) is exponentially stable.
 \end{corol}
{\bf Proof.} We choose the following equations: \\
$x(n+1)-x(n) = -a(n)x(g(n))$,
\\
$x(n+1)-x(n) = -a(n)x(g(n)) - b(n)x(g(n))$,
\\
with a positive fundamental function to obtain 1) and 2), respectively.
\qed

Consider now an autonomous equation with two delays
\begin{equation}
\label{28}
x(n+1)-x(n)=-ax(n-g)-bx(n-h), ~~ ag\not = 0, ~~ bh\not = 0.
\end{equation}
We further apply the nonoscillation condition (\ref{autonomnonoscil}).
\begin{corol}
\label{corol6}
Suppose at least one of the following conditions holds:
  
1) $ 0<a\leq g^g(g+1)^{-(g+1)},~ |b|<a $;
 
%
2) $ 0<(a+b)\leq g^g(g+1)^{-(g+1)},~ |a(g-h)|<1 $;


Then equation (\ref{28}) is exponentially stable.
 \end{corol}

Consider a high order autonomous difference equation.
\begin{equation}
\label{28a}
x(n+1)-x(n)=-\sum_{l=1}^m a_lx(n-l).
\end{equation} 
\begin{corol}
\label{corol11}
Suppose there exists $k\geq 1$ such that
$$
0<\sum_{l=1}^k a_l\leq \frac{k^k}{(k+1)^{(k+1)}},~~
\sum_{l=k+1}^m |a_l|<\sum_{l=1}^k a_l.
$$
Then equation (\ref{28a}) is exponentially stable.
\end{corol}

\section{Discussion and Examples}

The present work continues our previous publication \cite{JMAA2005}
where a first order exponentially stable model
$x(n+1)-x(n)= -b(n)x(n)$, $0<\varepsilon \leq b(n)\leq \gamma <1$
was used as a comparison equation. The results of the present paper
extend and improve most of theorems obtained in \cite{JMAA2005}.

Let us note that the approach using Bohl-Perron Theorem is similar to the 
method developed in \cite{Pituk2003} where stability is deduced based on
the fact that some linear exponentially stable equation is {\em close}
to the considered equation. Unlike the present paper,  \cite{Pituk2003} 
considers nonlinear perturbations of stable linear equations as well.
The main result (Theorem 2) of \cite{Pituk2003} is the following one.

{\em Suppose that the fundamental function of (\ref{1}) satisfies
\begin{equation}
\label{32a}
\sum_{j=0}^{n-1} |X(n,j+1)| \leq L, ~~n=n_0,n_0+1, \cdots ~.
\end{equation} 
Then the nonlinear equation
$$
x(n+1)-x(n)=-\sum_{k=1}^m a_k(n)x(h_k(n)) + F(n,x(n),x(n-1), \cdots, x(n-l))
$$
is globally asymptotically stable if in addition
${\displaystyle |F(n,x_0,x_1, \cdots , x_l)| 
\leq q \max_{0 \leq i \leq l} |x_i|}$, and $q<L^{-1}$.}

Instead of inequality (\ref{32a}) in this paper we apply exponential estimation (\ref{star}).
Generally, (\ref{star}) implies (\ref{32a}), however for bounded delays and coefficients
inequalities (\ref{star}) and (\ref{32a}) are equivalent. Indeed,
solution representation (\ref{3add}) and inequality (\ref{32a}) imply that 
for any bounded right hand side $|f(n)| \leq M$ the solution of the 
problem (\ref{1gen}),(\ref{2}) with the zero initial conditions
($\varphi(n) \equiv 0$, $x(n_0)=0$) is bounded: $|x(n)| \leq LM$.
Thus by the Bohl-Perron theorem (see Theorem 2 in \cite{JMAA2005}) the
fundamental function satisfies (\ref{star}), see also
Lemma 3 in \cite{Pituk2003} which claims exponential decay of solutions for autonomous equation 
(\ref{54}) if (\ref{32a}) holds and coefficients $a_k$ are positive. 


Let us discuss some stability tests for equation (\ref{1}).

We start with the following result \cite{ErbeXiaYu,MalKul,Yu,ZhTianChuan}.

{\em If $m=1$, 
$\sum_{n=0}^{\infty} a(n) = \infty$, $n-h(n)\leq k$ , $a(n)\geq 0$,
and
\begin{equation}
\label{32cond}
\sum_{i=h(n)}^{n} a(i)<\frac{3}{2}+\frac{1}{2k+2},
\end{equation}
then equation (\ref{1}) is asymptotically stable.
This result is also true for general equation (\ref{1}) ($m>1$), 
where $a_l(n)\geq 0$, $a(n)=\sum_{l=1}^m a_l(n), h(n)=\max h_l(n)$.}

Equation (\ref{1}) with positive constant coefficients is asymptotically stable
if \cite{GyoriHartung}
\begin{equation}
\label{32acond}
\sum_{l=1}^m a_l
~\limsup_{n\rightarrow\infty}
(n-h_l(n))<1+\frac{1}{e}-\sum_{l=1}^m a_l.
\end{equation}
Stability tests (\ref{32cond}) and  (\ref{32acond}) are obtained 
for equations with positive coefficients.
In the present paper we consider coefficients of arbitrary signs.
The next interesting feature of the results obtained here is that some of the
delays can be arbitrarily large 
(see for example, Parts 1 and 2 of Corollary \ref{corol5}). 

\begin{example}
By Corollary \ref{corol5}, Part 1, the following two equations
\begin{equation}
\label{51add}
x(n+1)-x(n)=-(0.2+0.05 \sin n)x(n-1)-0.1|\cos n| x(n-20),
\end{equation}
\begin{equation}
\label{52add}
x(n+1)-x(n)=-[0.12+0.1(-1)^n] x(n-2)-[0.1+0.11(-1)^n] x(n-14)
\end{equation}
are exponentially stable.

Two previous results 
of \cite{ErbeXiaYu,MalKul,Yu,ZhTianChuan} and \cite{GyoriHartung} fail
to establish exponential stability for equation (\ref{51add}) 
with positive coefficients, as well as all parts of
Corollary 3.10 in \cite{BBL} cannot be applied to equation (\ref{52add}) 
with an oscillating coefficient.
None of the inequalities in Corollary 8 of
\cite{JMAA2005} can be applied to establish stability of (\ref{52add}).

We note that (\ref{51add}),(\ref{52add}) are special cases
of equation with one nondelay term and two delay terms considered in   
\cite{JMAA2005}, however none of the inequalities in Corollary 8 of   
\cite{JMAA2005} can be applied to establish stability of (\ref{52add}).
Let us also note that for (\ref{52add}) we have
$$\sum_{k=g(n)}^{n-1} |a(k)|+\sum_{j=h(n)}^{n-1} |b(j)|=
\sum_{k=n-2}^{n-1} [0.12+0.1(-1)^k]+\sum_{j=n-14}^{n-1} 
\left| 0.1+0.11(-1)^k \right| =1.78>\frac{\pi}{2},
$$
where \cite{n21} $\pi/2$ is the best possible constant \cite{MalKul,TT} in
$\sum_{k=1}^m k a_k <\pi/2$ which implies exponential stability of (\ref{28a}).
\end{example}

\begin{example}
Consider equation (\ref{27}) with two variable coefficients and delays,
where
\begin{equation}
\label{eqexadd2}
a(n)=\left\{ \begin{array}{ll} -0.12,  & \mbox{if~} n \mbox{~is even}, \\
-0.05,  & \mbox{if~} n \mbox{~is odd}, \end{array} \right. ~
b(n)=\left\{ \begin{array}{ll}
0.17, & \mbox{if~} n \mbox{~is even}, \\
0.08, & \mbox{if~} n \mbox{~is odd}, \end{array} \right.
\end{equation}
\begin{equation}
\label{eqexadd2a}
g(n)=\left\{ \begin{array}{ll} n-3, & \mbox{if~} n \mbox{~is even}, \\
n-5, & \mbox{if~} n \mbox{~is odd} \end{array} \right. ~
h(n)=\left\{ \begin{array}{ll} n-4, & \mbox{if~} n \mbox{~is even}, \\   
n-8, & \mbox{if~} n \mbox{~is odd.} \end{array} \right.
\end{equation}  
Then in (\ref{27}) the sum $a(n)+b(n)$ is either 0.05 or 0.03 which is less
than 0.5, $\sum_{k=n-5}^{n-1} [a(k)+b(k)]\leq 3 \cdot 0.05+2 \cdot 0.03 =
0.21<1/4$
for both odd and even $n$ and for $\gamma=0.95<1$ we have    
$$
|a(n)| \sum_{k=\min\{h(n),g(n)\}}^{\max\{h(n),g(n)\}-1} [|a(k)|+|b(k)|]
= \left\{ \begin{array}{ll} 0.12 \cdot 0.29=0.0348  < 0.05 \gamma, & n
\mbox{~is even}, \\
0.05 [0.13 \cdot 2+0.29]=0.0275 < 0.03 \gamma, & n \mbox{~is odd},
\end{array} \right.
$$
thus by Corollary \ref{corol5}, Part 3, equation (\ref{27}) is
exponentially stable. None of the inequalities in Corollary 8 of   
\cite{JMAA2005} can be applied to establish stability of (\ref{27}).
We note that it would be harder to treat this example if the 
equation were written as high order equations with constant
delays and variable coefficients.
\end{example}

A  number of papers \cite{JMAA2005,GLV,KL,KP,L,LIF,Pituk2003} are devoted to 
stability 
tests for equations 
with positive and negative coefficients and, more generally, for equations with 
oscillating coefficients. Paper \cite{Pituk2003} extends 
earlier results of \cite{GLV}.
In particular, for the linear autonomous equation
\begin{equation}
\label{2aut}
x(n+1)-x(n)=qx(n-m)-px(n-k),~ p>0, q>0, m\geq 1, k\geq 1.
\end{equation}
the following result was obtained in \cite{GLV}.

{\em Suppose $p\frac{(k+1)^{(k+1)}}{k^k}\leq 1$. Then equation(\ref{2aut}) 
is exponentially stable if and only if $p>q$. }

Condition 1) of Corollary \ref{corol6} is close to this result.
It gives the same sufficient stability test for $q$ of an 
arbitrary sign but does not involve the necessity part.

The paper \cite{L} contains a nice review on stability results obtained
for equations with oscillating coefficients. 
The results of \cite{L} generalized the following stability test
obtained in \cite{KL} for equation (\ref{2aut}):

{\em If $kp<1, p\frac{1-kp}{1+kp}>q$ 
then (\ref{2aut}) is asymptotically stable. }

By condition 2) of Corollary \ref{corol6} equation (\ref{2aut}) is 
asymptotically stable if
$$p-q<\frac{k^k}{(k+1)^{(k+1)}},~~ |p(k-m)|<1.$$
It is easy to see that these two tests are independent.

Let us discuss sharpness of conditions of Theorem 1
for exponential stability
of (\ref{1}), assuming the fundamental function is positive;
in particular, we demonstrate sharpness of condition (\ref{16add}).

\begin{example}
The equation
$$
x(n+1)-x(n)=-3^{-n-1}x(n),~~n\geq n_0\geq 0,
$$
has a positive fundamental function
and any solution can be presented as
\\
${\displaystyle x(n) = x(n_0) \prod_{k=n_0}^{n-1} (1-3^{-k-1}) }$, thus
for $n_0=0$
$$X(n,0)=\prod_{k=0}^{n-1} (1-3^{-k-1}) > 1 - \sum_{k=0}^{n-1} 3^{-k-1} >
1 - \sum_{k=0}^{\infty} 3^{-k-1}= \frac{1}{2},$$
i.e., the equation is neither asymptotically nor exponentially stable.
\end{example}

Let us demonstrate that the facts that the sum of
coefficients ${\displaystyle \sum_{l=1}^m a_l(n)}$  in (\ref{1}) is
positive, exceeds a positive number and that the
fundamental function is positive
do not imply stability, in the case when coefficients have different
signs.
\vspace{2mm}

\begin{example}
Consider the difference equation
\begin{equation}
\label{eqexample}
x(n+1)-x(n)=-2.2x(n-1)+2x(n).
\end{equation}
Here $2.2-2=0.2>0$, so the sum of coefficients exceeds a certain positive
number. Let us prove that the
fundamental function is positive and the solution is unbounded.
Really, for the fundamental function we have $X(0,0)=1, X(1,0)=3$. Denote
$x(n)=X(n,0)$, notice that $x(1)>1.5 x(0)$ and prove $x(n)>1.5x(n-1)>0$ by
induction.
Really, $x(n)>1.5x(n-1)>0$ implies $x(n-1)<2x(n)/3$, and
for any $x(n-1)>0$ we have 
$$x(n+1)=3x(n)-2.2x(n-1)>3x(n)-4.4x(n)/3=\frac{4.6}{3}x(n)>1.5 x(n),$$
thus $X(n,0)$ is positive and unbounded. The equation is autonomous, so
the same is true for $X(n,k)$. Since $X(n,0)$ is unbounded, then
(\ref{eqexample}) is not stable.
\end{example}

Finally, let us formulate some open problems.

\begin{enumerate}
\item
Under which conditions will exponential stability of (\ref{1}) 
imply exponential stability of the equation with the same coefficients and 
smaller delays:
$$
x(n+1)-x(n) = - \sum_{l=1}^m a_l(n) x(g_l(n)), ~n \geq n_0,~h_l(n) \leq 
g_l(n) \leq n ~?$$
\item 
Prove or disprove:\\
If in Theorem \ref{dominate} condition (\ref{55}) is substituted by
$$\sum_{l \not\in I} |a_l(n)| \leq \alpha_n \sum_{l \in I} a_l(n),
~~ \prod_{n=1}^{\infty} \alpha_n \leq 1$$
and all other assumptions hold, then (\ref{1}) is stable. If in addition
$$ \sum_{n=1}^{\infty} (1-\alpha_n)=\infty, $$
then (\ref{1}) is asymptotically stable.
\item
Consider the problem of the exponential stability of (\ref{1}) when
(a1)-(a2) are substituted with one of two more general conditions:
\begin{description}
\item[a)] ${\displaystyle \lim_{n \to \infty} h_l(n) = \infty}$ and there 
exists
$M>0$ such that
$$ \sum_{j=h_l(n)}^n a_l(j) < M,~ \mbox{ for any }~ l=1, \cdots , m; $$

\item[b)] delays are infinite but coefficients decay
exponentially with memory, i.e., there exist positive numbers $M$ and
$\lambda<1$ such that ${\displaystyle |a_l(n)| \leq M \lambda^{n-h_l(n)} 
}$.
\end{description}

Let us note that Bohl-Perron type result in case b) was obtained in
\cite{melbourne}, Theorem 4.7.
\item
Example 5 demonstrates that for equations with positive and negative
coefficients and a positive fundamental function
inequality (\ref{16add}) does not imply exponential stability. Is it
possible to find such conditions on delays and coefficients of different  
signs that (\ref{16add}) would imply exponential stability? For 
instance, prove or disprove the following conjecture.

{\bf Conjecture.} Suppose the following conditions
$$a(n) \geq b(n) \geq 0, ~~h(n) \leq g(n) \leq n,
~~ \limsup_{n \to \infty} b(n)[g(n)-h(n)]<1 $$
are satisfied for the equation
\begin{equation}
\label{flower}   
x(n+1)-x(n)=-a(n)x(h(n))+ b(n) x(g(n)).
\end{equation}
If the fundamental function of (\ref{flower}) is positive and
$$ \liminf_{n \to \infty} [a(n)-b(n)] >0,$$
then (\ref{flower}) is exponentially stable.
\\
\\
Let us remark that conditions when the fundamental function of
(\ref{flower}) is positive were obtained in \cite{oksana}.
\end{enumerate}

\vspace{3mm}

\centerline{\large\bf Acknowledgement}
\vspace{3mm}

The authors are grateful to the referee for valuable comments and remarks.


\begin{thebibliography}{99}


\bibitem{FDE2004}
L. Berezansky and E. Braverman, On Bohl-Perron type theorems for linear 
difference equations, {\em Funct. Differ. Equ.} {\bf 11}  (2004),  no. 
1-2, 19--28.
 
\bibitem{JMAA2005}
L. Berezansky and E. Braverman, 
On exponential dichotomy, Bohl-Perron type theorems and stability of
difference equations, {\em J. Math. Anal. Appl.} {\bf 304} (2005), 
511--530.

\bibitem{ADSA}
L. Berezansky and E. Braverman,
On existence of positive solutions for linear difference equations with   
several delays, 
{\em  Adv. Dyn. Syst. Appl.} {\bf 1}  (2006),  no. 1, 29--47.

\bibitem{BBL}
L. Berezansky, E. Braverman and E. Liz,  Sufficient conditions for the
global stability of nonautonomous higher order difference equations, {\em
J. Difference Equ. Appl.} {\bf 11} (2005), no. 9, 785--798.

\bibitem{melbourne}
L. Berezansky and E. Braverman,
On exponential dichotomy for linear difference equations with bounded and
unbounded delay, in {\em Proceedings of the Conference on Differential
and Difference Equations \& Applications} Melbourne, Florida, August 1-5, 
2005. Edited by
Ravi P. Agarwal and Kanishka Perera, Hindawi Publishing Corporation, 2006, 
169--178.

\bibitem{oksana}
L. Berezansky, E. Braverman and O. Kravets, Nonoscillation of linear delay
difference equations with positive and negative coefficients, 
{\em J. Difference Equ. Appl.} {\bf 14}  (2008),  no. 5, 495--511.

\bibitem{E}
S. Elaydi, Periodicity and stability of linear Volterra difference 
systems, {\em J. Math. Anal. Appl.} {\bf 181} (1994), no. 2, 
483--492.

\bibitem{n20}
S. Elaydi and S. Murakami, Uniform asymptotic stability in linear 
Volterra difference equations, {\em J. Difference Equ. Appl.} {\bf 3} 
(1998),  no. 3-4, 203--218. 


\bibitem{ErbeXiaYu}
L.\,H. Erbe,  H. Xia and J.\,S. Yu,
Global stability of a linear nonautonomous delay
difference equation, {\em J. Difference Equ. Appl.}
{\bf  1} (1995), no. 2, 151--161.

\bibitem{GyoriHartung}
I. Gy\"{o}ri and F. Hartung, Stability in delay perturbed differential and
difference equations, {\em Fields Inst. Commun.} {\bf 29}
(2001), 181--194.

\bibitem{GL}
I. Gy\"{o}ri and G. Ladas, 
Oscillation theory of delay differential equations.
Clarendon Press, Oxford, 1991.

\bibitem{GLV}
I. Gy\"{o}ri, G. Ladas and P.\,N. Vlahos,
Global attractivity in a delay difference equation,
{\em Nonlinear Anal. TMA} {\bf 17} (1991), no. 5, 473--479.


\bibitem{n21}
M.\,M. Kipnis and D.\,A. Komissarova,
A note of explicit stability conditions of autonomous higher order
difference equation, {\em J. Difference Equ. Appl.} {\bf 13} (2007), 
no. 5, 457--461.

\bibitem{KL} 
V.\,L. Koci\'{c} and G. Ladas, Global Behavior of Nonlinear Difference Equations of 
Higher Order with Applications,
Math. Appl., vol. 256, Kluwer Academic, Dordrecht, 1993.

\bibitem{KP} 
U. Krause and M. Pituk, Boundedness and stability for higher order 
difference equations, 
{\em J. Differ. Equations Appl.} {\bf 10} (2004), 343--356.

\bibitem{L} 
E. Liz, On explicit conditions for the asymptotic stability
of linear higher order difference equations, 
{\em J. Math. Anal. Appl.}  {\bf 303} (2005), 492--498.

\bibitem{LIF} 
E. Liz, A. Ivanov, J.B. Ferreiro, 
Discrete Halanay-type inequalities and applications, 
{\em Nonlinear Anal.} {\bf 55}
(2003), 669--678.

\bibitem{n8}
E. Liz and M. Pituk,  Asymptotic estimates and exponential 
stability for higher-order monotone difference equations, {\em  Adv. 
Difference Equ.}  2005,  no. 1, 41--55. 

\bibitem{MalKul} 
V.\,V. Malygina and A.\,Y. Kulikov,
On precision of constants in some theorems on stability of
difference equations,
{\em Func. Differ. Equ.} {\bf 15} (2008), no. 3-4, 239--249.




\bibitem{Pituk2003}
M. Pituk, Global asymptotic stability in a perturbed higher order 
linear difference equation, {\em Comput. Math. Appl.} {\bf 45} (2003), 
1195--1202.

\bibitem{n22}
X.\,H. Tang and Z. Jiang, Asymptotic behavior of Volterra 
difference equation, {\em J. 
Difference Equ. Appl.} {\bf 13} (2007), no. 1, 25--40.

\bibitem{TT}
V. Tkachenko and S. Trofimchuk, A global attractivity criterion 
for nonlinear non-autonomous difference equations, {\em J. Math. Anal. 
Appl.} {\bf 322}  (2006),  no. 2, 901--912.

\bibitem{Yu}
J.\,S. Yu, Asymptotic stability for a linear difference equation with 
variable delay,
Advances in difference equations, II, {\em Comput. Math. Appl.} {\bf 36}
(1998), no. 10-12, 203--210.  

\bibitem{ZhTianChuan}
B.\,G. Zhang, C.\,J. Tian and P.\,J.\,Y. Wong,
Global attractivity of difference equations with variable delay,
{\em Dynam. Contin. Discrete Impuls. Systems} {\bf 6} (1999), no. 3,
307--317.


\end{thebibliography}
\end{document}